\documentclass{amsart}
\usepackage[OT2,T1]{fontenc}
\usepackage{lmodern}
\usepackage{amsmath,amssymb,amsfonts}
\usepackage[pdfborder={0 0 0}]{hyperref}
\usepackage[ps,dvips,arrow,matrix,tips,line]{xy}
\entrymodifiers={+!!<0pt,\the\fontdimen22\textfont2>}
\SelectTips{cm}{10}
\input cyracc.def
\DeclareFontFamily{U}{russian}{}
\DeclareFontShape{U}{russian}{m}{n}
        { <5><6> wncyr5
        <7><8><9> wncyr7
        <10><10.95><12><14.4><17.28><20.74><24.88> wncyr10 }{}
\DeclareSymbolFont{Russian}{U}{russian}{m}{n}
\DeclareSymbolFontAlphabet{\mathcyr}{Russian}
\makeatletter
\let\@math@cyr\mathcyr
\renewcommand{\mathcyr}[1]{\@math@cyr{\cyracc #1}}
\makeatother
\newcommand{\Ba}{{\mathcyr{B}}}
\newtheorem{thm}{Theorem}[section]
\newtheorem{lem}[thm]{Lemma}
\newtheorem*{lem*}{Lemma}
\newtheorem{cor}[thm]{Corollary}

\theoremstyle{definition}

\theoremstyle{remark}
\newtheorem{rmk}[thm]{Remark}
\newtheorem{rmks}[thm]{Remarks}
\newtheorem{question}[thm]{Question}
\numberwithin{equation}{section}
\DeclareMathOperator{\inv}{inv}
{\setbox0\hbox{$ $}}\fontdimen16\textfont2=\fontdimen17\textfont2
\renewcommand{\emptyset}{\varnothing}
\makeatletter
\def\myrightarrow{{\setbox\z@\hbox{$\rightarrow$}\dimen0\ht\z@\multiply\dimen0 6\divide\dimen0 10\ht\z@\dimen0\box\z@}}
\def\myrightarrowfill@{\arrowfill@\relbar\relbar\myrightarrow}
\newcommand{\myxrightarrow}[2][]{\ext@arrow 0359\myrightarrowfill@{#1}{#2}}
\makeatother
\newcommand{\isoto}{\myxrightarrow{\,\sim\,}}
\def\tilde{\widetilde}
\newcommand{\bigset}[2]{\left\{ #1 \, ; \, #2 \right\}}
\newcommand{\mylangle}{\xy(0,0),(1.17,1.485)**[|(1.025)]@{-},(0,0),(1.17,-1.485)**[|(1.025)]@{-}\endxy\mkern2mu}
\newcommand{\myrangle}{\mkern2mu\xy(0,0),(-1.17,1.485)**[|(1.025)]@{-},(0,0),(-1.17,-1.485)**[|(1.025)]@{-}\endxy}
\newcommand{\Br}{{\mathrm{Br}}}
\newcommand{\Gm}{\mathbf{G}_\mathrm{m}}
\newcommand{\Zhat}{{\hat \Z}}
\newcommand{\kbar}{{\mkern1mu\overline{\mkern-1mu{}k\mkern-1mu}\mkern1mu}}
\newcommand{\Kbar}{{\mkern2mu\overline{\mkern-2mu{}K}}}
\newcommand{\kbarstar}{\vphantom{\kbar}\smash{\smash{\kbar}^\star}\vphantom{\kbar}}
\newcommand{\Kbarstar}{\vphantom{\kbar}\smash{\smash{\Kbar}^\star}\vphantom{\kbar}}
\newcommand{\kXbar}{\overline{k(X)}}
\newcommand{\Alb}{\mathrm{Alb}}
\newcommand{\mmu}{{\boldsymbol{\mu}}}
\newcommand{\Pic}{{\rm Pic}}
\newcommand{\Div}{{\rm Div}}
\newcommand{\et}{{\text{\'et}}}
\newcommand{\Galp}[1]{{\rm Gal}(#1)}
\newcommand{\Galpp}[1]{{\rm Gal}\mkern-2.5mu\left(#1\mkern-1mu\right)}
\newcommand{\C}{{\mathbf C}}
\newcommand{\Q}{{\mathbf Q}}
\newcommand{\Z}{{\mathbf Z}}

\newcommand{\ab}{{\mathrm{ab}}}
\newcommand{\Gk}{{\vphantom{G_{k(X)}}\smash{G_{\vphantom{k(X)}\mkern-1mu{}k}^{\vphantom{[ab]}}}\vphantom{G_{k(X)}}}}
\newcommand{\GK}{{\vphantom{G_{K(X)}}\smash{G_{\vphantom{K(X)}\mkern-1mu{}K}^{\vphantom{[ab]}}}\vphantom{G_{K(X)}}}}

\newcommand{\GkX}{{\vphantom{G_{k(X)}}\smash{G_{k(X)}^{\vphantom{[ab]}}}\vphantom{G_{k(X)}}}}
\newcommand{\Ggeom}{\vphantom{G_{\smash{\kbar}(X)}}\smash{G^{\vphantom{]}}_{\smash{\kbar}(X)}}\vphantom{G_{\smash{\kbar}(X)}}}
\newcommand{\Gabgeom}{G^{{\mkern1mu\mathrm{ab}\vphantom{]}}}_{\smash{\kbar}(X)}}
\newcommand{\Gabgeomprimetop}{G^{{\mkern1mu\mathrm{ab},p'\vphantom{]}}}_{\smash{\kbar}(X)}}
\newcommand{\Gabgeomp}{G^{{\mkern1mu\mathrm{ab},p\vphantom{]}}}_{\smash{\kbar}(X)}}
\newcommand{\Gabgeomell}{G^{{\mkern1mu\mathrm{ab},\ell\vphantom{]}}}_{\smash{\kbar}(X)}}
\newcommand{\Gab}{G^{\mathrm{[ab]}}_{k(X)}}
\newcommand{\GabK}{G^{\mathrm{[ab]}}_{K(X)}}
\newcommand{\ob}{\mathrm{ob}}
\newcommand{\Ext}{{\mathrm{Ext}}}
\newcommand{\Ker}{{\mathrm{Ker}}}
\newcommand{\Hom}{{\mathrm{Hom}}}
\newcommand{\Cores}{{\mathrm{Cores}}}
\title{On abelian birational sections}
\author{H\'el\`ene Esnault}
\address{Universit\"at Duisburg--Essen, Mathematik, 45117 Essen, Germany}
\email{esnault@uni-due.de}
\author{Olivier Wittenberg}
\address{D\'epartement de math\'ematiques et applications, \'Ecole normale sup\'erieure, 45~rue d'Ulm, 75320 Paris Cedex 05, France}
\email{wittenberg@dma.ens.fr}
\thanks{Supported in part by the DFG Leibniz Preis, the SFB/TR 45 and the ERC/Advanced Grant~226257}
\date{February 9, 2009; revised on November 27, 2009}
\subjclass[2000]{Primary 14G32; Secondary 14C25, 14G25, 14G20}

\begin{document}
\vspace*{-.9mm}

\begin{abstract}
For a smooth and geometrically irreducible variety~$X$ over a field~$k$,
the quotient of the absolute Galois group $G_{k(X)}$ by
the commutator subgroup of $G_{\bar k(X)}$ projects onto~$G_k$.
We investigate the sections of this projection.  We show that such sections
correspond to ``infinite divisions'' of the elementary obstruction of
Colliot-Th\'el\`ene and Sansuc.  If~$k$ is a number field and the
Tate--Shafarevich group of the Picard variety of~$X$ is finite, then such
sections exist if and only if the elementary obstruction vanishes.  For
curves this condition also amounts to the existence of divisors of
degree~$1$.  Finally we show that the vanishing of the elementary
obstruction is not preserved by extensions of scalars.
\end{abstract}

\maketitle

\section{Introduction}

In his seminal letter to Faltings, Grothendieck~\cite{gtof} formulated the conjecture, now known as the section conjecture, that
for any smooth and proper curve~$X$ of genus~$\geq 2$
defined over a field~$k$ finitely generated over~$\Q$, sections of the natural exact sequence of fundamental groups
\begin{equation}
\label{introsepi}
\xymatrix{
1 \ar[r] & \pi_1(X \otimes_k \kbar) \ar[r] & \pi_1(X) \ar[r] & \Galp{\kbar/k} \ar[r] & 1
}
\end{equation}
are, up to conjugation by $\pi_1(X \otimes_k \kbar)$, in one-to-one correspondence with rational points of~$X$.
Accordingly, the exact sequence of profinite groups~(\ref{introsepi}) should split if and only if $X(k) \neq \emptyset$;
this corollary is in fact equivalent to the full statement of the section conjecture (see~\cite[Lemma~1.7]{koenigsmann}).

As was first emphasised by Koenigsmann~\cite{koenigsmann}, a birational variant of the section conjecture may be stated as follows:
for~$k$ and~$X$ as above,
if we denote by~$G_K$ the absolute Galois group of a field~$K$,
the exact sequence of profinite groups
\begin{equation}
\label{introsegnonab}
\xymatrix{
1 \ar[r] & \Ggeom \ar[r] & \GkX \ar[r] & \Gk \ar[r] & 1
}
\end{equation}
should split if and only if~$X(k)\neq\emptyset$.
That~(\ref{introsegnonab}) must split if $X(k)\neq\emptyset$ is a general remark of Deligne~\cite[Section~15]{delignegroupefond}.
The converse implication would follow from the original section conjecture, since splittings of~(\ref{introsegnonab})
induce splittings of~(\ref{introsepi}).  To this day even the birational section conjecture is widely open.

The goal of this note is to investigate the exact sequence of profinite groups
\begin{equation}
\label{introsegab}
\xymatrix{
1 \ar[r] & \Gabgeom \ar[r] & \Gab \ar[r] & \Gk \ar[r] & 1
}
\end{equation}
\vskip-.5em{}\noindent{}obtained by pushing out~(\ref{introsegnonab}) by the abelianisation map $\Ggeom \rightarrow \Gabgeom$.
We~refer to splittings of~(\ref{introsegab}) as ``abelian birational sections''.
Abelian birational sections exist whenever~$X$ possesses a divisor of degree~$1$ defined over~$k$;
the question therefore arises whether the converse might hold, and if so, in what generality.

In \textsection\ref{sectionbirat} we prove (Theorem~\ref{thmab})
that if~$k$ is a number field and~$X$ is a smooth and proper curve over~$k$ whose Jacobian has finite Tate--Shafarevich group,
then~(\ref{introsegab}) splits if and only if~$X$ possesses a divisor of degree~$1$.
This is an abelian variant of the birational section conjecture over number fields.
It should be noted that starting from the existence of an abelian birational section, one cannot hope to prove more than the existence of a divisor of degree~$1$,
since there are smooth and proper curves over~$\Q$ admitting a divisor of degree~$1$ but no rational point.

We then proceed in~\textsection\ref{sectob} to study abelian birational sections for an arbitrary smooth variety~$X$ over an arbitrary field~$k$.

In this situation, the existence of an abelian birational section is a necessary condition for the existence of a $0$\nobreakdash-cycle of degree~$1$ on~$X$.
Remarkably few necessary conditions for the existence of $0$\nobreakdash-cycles of degree~$1$ have been exhibited in such generality.
Another one is the vanishing of the elementary obstruction of Colliot-Th\'el\`ene and Sansuc (see~\cite{ctsandesc2}, \cite{skobook},
\cite{boctsk}, \cite{w}), which~takes the form of a class $\ob(X)$ in $\Ext^1_{G_k}(\kbar(X)^\star/\kbarstar,\kbarstar)$.
In~\textsection\ref{sectob} we determine the precise relation between these two obstructions.
Namely the main result of~\textsection\ref{sectob} (Theorem~\ref{thmcomp}) states
that for any field~$k$ of characteristic~$0$ and any variety~$X$ over~$k$, the existence of an abelian birational section is equivalent to~$\ob(X)$
belonging to the maximal divisible subgroup of $\Ext^1_{G_k}(\kbar(X)^\star/\kbarstar,\kbarstar)$
(see Remark~\ref{imperfect} for the situation in positive characteristic).
In particular~(\ref{introsegab}) splits whenever the elementary obstruction vanishes.
In the course of the proof of Theorem~\ref{thmcomp}, we establish a one-to-one correspondence
between abelian birational sections, up to conjugation by~$\smash{\Gabgeom}$, and isomorphism classes of discrete $G_k$\nobreakdash-modules which are simultaneously extensions
of $\left(\kbar(X)^\star/\kbarstar\right) \otimes_\Z \Q$ by~$\kbarstar$ and of $\left(\kbar(X)^\star/\kbarstar\right) \otimes_\Z \Q/\Z$
by $\kbar(X)^\star$ (Lemma~\ref{lemsecconj}).

Now assume that~$k$ is a number field.
From Theorem~\ref{thmcomp} it follows that the statement of Theorem~\ref{thmab} does not hold for varieties of arbitrary dimension.
Indeed it is well known that there are varieties over number fields (for instance, rational surfaces over~$\Q$)
which do not possess $0$\nobreakdash-cycles of degree~$1$
although the elementary obstruction vanishes
(in which case, by Theorem~\ref{thmcomp}, abelian birational sections do exist).
Yet Grothendieck writes in~\cite[p.~3]{gtof} that any sufficiently small nonempty open subset of a smooth variety should be ``anabelian'',
which, according to the general statement of the section conjecture given in~\cite{gtof}, implies that
the birational section conjecture should hold for smooth varieties of arbitrary dimension.
The correct higher-dimensional generalisation of Theorem~\ref{thmab} is given by Theorem~\ref{thmobobdiv} below,
which refines, for number fields, the statement of Theorem~\ref{thmcomp}: if~$k$ is a number field and~$X$ is
a smooth and proper variety over~$k$ whose Picard variety has finite Tate--Shafarevich group, then the existence of an abelian birational section
is equivalent to the vanishing of the elementary obstruction.

Finally, by considering sections of the exact sequence~(\ref{introsepi}) and their cycle classes in \'etale cohomology (see~\cite{ew}),
we give in Theorem~\ref{thmcex}
 an example of a field~$k$ of characteristic~$0$, a field extension~$K/k$, and a curve~$X$ over~$k$ such that $\ob(X)=0$ but $\ob(X \otimes_k K)\neq 0$.
This settles a question about the elementary obstruction which was raised by Borovoi, Colliot-Th\'el\`ene and Skorobogatov~\cite{boctsk}.

\section{Abelian birational sections for curves over number fields}
\label{sectionbirat}

Let~$k$ be a field and~$\kbar$ be a separable closure of~$k$.  Let~$X$ be a smooth, proper, geometrically connected curve over~$k$.
We denote by $k(X)$ the function field of~$X$, by $\kXbar$ a separable closure of~$\kbar(X)=k(X) \otimes_k \kbar$
and by $\kbar(X)^\ab$ the maximal abelian subextension of $\kXbar/\kbar(X)$.
The Galois group $\Gab=\Galpp{\kbar(X)^\ab/k(X)}$ fits into an exact sequence of profinite groups
\begin{equation}
\label{seg}
\xymatrix{1 \ar[r] & \Gabgeom \ar[r] & \Gab \ar[r] & \Gk \ar[r] & 1\text{,}}
\end{equation}
where $\Gk=\Galp{\kbar/k}$ and $\Gabgeom=\Galpp{\kbar(X)^\ab/\kbar(X)}$.
Equivalently, the exact sequence~(\ref{seg}) is the one obtained by pushing out
\begin{equation}
\label{segnonab}
\xymatrix{1 \ar[r] & \Ggeom \ar[r] & \GkX \ar[r] & \Gk \ar[r] & 1\text{,}}
\end{equation}
where $\Ggeom=\Galpp{\kXbar/\kbar(X)}$ and $\GkX=\Galpp{\kXbar/k(X)}$, by the abelianisation
map $\Ggeom\rightarrow \Gabgeom$.

According to Deligne~\cite[Section~15]{delignegroupefond} (see also Lemma~\ref{deligne} below), the exact sequence~(\ref{segnonab}) splits if $X(k)\neq\emptyset$.
From this it follows, by a standard restriction-corestriction argument taking place in the continuous
Galois cohomology group $H^2(k,\Gabgeom)$, that~(\ref{seg}) splits as soon as~$X$ possesses a degree~$1$ divisor.
In this section we prove that over number fields, the converse holds, under the assumption that the Tate--Shafarevich group of the Jacobian of~$X$ is finite:

\begin{thm}
\label{thmab}
Let~$X$ be a smooth proper geometrically connected curve over a number field~$k$.  Assume that the Tate--Shafarevich group of the Jacobian of~$X$ is finite.
Then the exact sequence of profinite groups~(\ref{seg}) splits if and only if there exists a divisor of degree~$1$ on~$X$.
\end{thm}

Recall that the birational section conjecture (a consequence of Grothendieck's
section conjecture in anabelian geometry)
asserts
 over number fields that
for~$X$ and~$k$ as above, the exact sequence of profinite
groups~(\ref{segnonab}) splits if and only if there exists a rational point on~$X$.
(This conjecture is sometimes referred to as the ``weak birational section conjecture'', although, as a consequence of Faltings's theorem,
it is in fact equivalent to the full birational section conjecture; see~\cite[Lemma~1.7]{koenigsmann}.)
The statement of Theorem~\ref{thmab} may thus be considered as an abelian variant of the birational section conjecture over number fields ---
indeed the set of divisors of degree~$1$ on~$X$ is ``abelian'' in the sense that it is canonically a principal homogeneous
space under an abelian group, namely, the group of divisors of degree~$0$.

\begin{proof}[Proof of Theorem~\ref{thmab}]
The discrete $\Gab$\nobreakdash-module defined by
$$
M = \bigset{f \in \kbar(X)^{\ab \mkern1mu \star}}{\exists\mkern1mu n \geq 1, \mkern2mu{}f^n \in \kbar(X)^\star}
$$
naturally fits into the commutative diagram
\begin{equation}
\begin{aligned}
\label{diagrammem}
\xymatrix{
0 \ar[r] & \kbarstar \ar@{=}[d] \ar[r] & \kbar(X)^\star \ar[d]^i \ar[r] & \ar[d]\kbar(X)^\star/\kbarstar \ar[r] & 0 \\
0 \ar[r] & \kbarstar \ar[r] & M \ar[r]^(0.3)p & \left(\kbar(X)^\star/\kbarstar\right) \otimes_\Z \Q \ar[r] & 0 \rlap{\text{,}}
}
\end{aligned}
\end{equation}
where $i(f)=f$ and $p(f)=f^n \otimes \frac 1n$ for any~$n$ such that $f^n\in \kbar(X)^\star$.  
As Kummer extensions are abelian, the map~$p$ is onto.  Hence the rows of~(\ref{diagrammem}) are exact.

Suppose~(\ref{seg}) splits, and fix a section $s \colon \Gk \rightarrow \Gab$.  By letting~$G_k$ act on~$M$ \emph{via}~$s$,
we can view~(\ref{diagrammem}) as a commutative diagram of discrete $G_k$\nobreakdash-modules with exact rows.
Since the Galois cohomology of a discrete $G_k$\nobreakdash-module which at the same time is a $\Q$\nobreakdash-vector space vanishes in positive degrees,
the bottom row of~(\ref{diagrammem}) induces in Galois cohomology an isomorphism $H^2(k,\kbarstar) \isoto H^2(k,M)$.
On composing the inverse of this isomorphism with the map $H^2(k,\kbar(X)^\star) \rightarrow H^2(k,M)$ induced by~$i$, we find
a retraction of the natural map $H^2(k,\kbarstar) \rightarrow H^2(k,\kbar(X)^\star)$.  Now $H^2(k,\kbarstar)=\Br(k)$
and, by Tsen's theorem, $H^2(k,\kbar(X)^\star)=\Br(k(X))$.  Moreover, the Brauer group of the curve~$X$ is a subgroup of $\Br(k(X))$
(see~\cite[Cor.~1.8]{grothbr2}).
Thus, by restriction, we obtain a retraction $r \colon \Br(X) \rightarrow \Br(k)$ of the natural map $\Br(k) \rightarrow \Br(X)$.

Let~$v$ be a place of~$k$.  Choose a place of~$\kbar$ above~$v$, denote by $k_v^h \subset \kbar$ the fixed field of
the corresponding decomposition subgroup $D_v \subset G_k$, and consider~(\ref{diagrammem}) as a diagram of
discrete $D_v$\nobreakdash-modules. 
As in the previous paragraph,
passing to cohomology
yields
a retraction $r_v \colon \Br(X \otimes_k k_v^h) \rightarrow \Br(k_v^h)$ of the
natural map $\Br(k_v^h) \rightarrow \Br(X \otimes_k k_v^h)$.

Let~$k_v$ denote the completion of~$k$ at~$v$.
As is well known, the Brauer groups of~$k_v$ and of~$k_v^h$ coincide: by local class field theory,
both are canonically isomorphic to $\Q/\Z$ if~$v$ is finite
(see~\cite[Ch.~I, Prop.~A.1]{milneadt}) or to $\Z/2\Z$ if~$v$ is real.

\begin{lem}
\label{lemmakvhkv}
For any place~$v$ of~$k$, the natural map $\Br(X \otimes_k k_v^h) \rightarrow \Br(X \otimes_k k_v)$ is an isomorphism.
\end{lem}

\begin{proof}
The injectivity of this map results from theorems of Greenberg (for finite~$v$) and of Artin (for real~$v$),
see \cite[p.~334]{boctsk}.  To prove its surjectivity, embed~$\kbar$ into an algebraic closure~$K$ of~$k_v$.
The two field extensions $K/k_v$ and $\kbar/k_v^h$ then have the same Galois group, namely~$D_v$.
Consider the Hochschild--Serre spectral sequences in \'etale cohomology
$H^p(D_v,H^q_\et(X \otimes_k \kbar,\Gm)) \Rightarrow H^{p+q}_\et(X \otimes_k k_v^h,\Gm)$
and $H^p(D_v,H^q_\et(X \otimes_k K,\Gm)) \Rightarrow H^{p+q}_\et(X \otimes_k k_v,\Gm)$.
Since the Brauer groups of $X \otimes_k \kbar$ and of $X \otimes_k K$ vanish (by Tsen's theorem), one obtains a commutative diagram with exact rows
\begin{equation}
\begin{aligned}
\label{diagbr}
\xymatrix{
\Br(k_v^h) \ar[d]^\wr \ar[r] &\Br(X \otimes_k k_v^h) \ar[d] \ar[r] & H^1(D_v, \Pic(X \otimes_k \kbar)) \ar[d] \ar[r] & H^3(D_v,\kbarstar)  \\
\Br(k_v) \ar[r] &\Br(X \otimes_k k_v) \ar[r] & H^1(D_v, \Pic(X \otimes_k K))\rlap{\text{.}}
}
\end{aligned}
\end{equation}
Let~$J$ denote the Jacobian of~$X$.
The natural map $\Pic(X \otimes_k \kbar) \to \Pic(X \otimes_k K)$ is injective and its cokernel identifies with $J(K)/J(\kbar)$.
As the inclusion $J(\kbar) \subset J(K)$ induces an isomorphism on torsion subgroups
and both~$J(\kbar)$ and~$J(K)$ are divisible groups, the quotient $J(K)/J(\kbar)$ is a $\Q$\nobreakdash-vector space.
It follows that the rightmost vertical map appearing in~(\ref{diagbr}) is onto.
Moreover, the group $H^3(D_v,\kbarstar)$ vanishes (by \cite[Ch.~II, \textsection5.3, Prop.~15]{serrecg} for finite~$v$ and by Hilbert's Theorem~90 for real~$v$)
and the leftmost vertical map is an isomorphism.
Hence the middle vertical map is onto as well.
\end{proof}

We resume the proof of Theorem~\ref{thmab}.
Thanks to Lemma~\ref{lemmakvhkv}, we may view~$r_v$ as a retraction of the natural map
$\Br(k_v)\rightarrow \Br(X\otimes_k k_v)$.
Every divisor of degree~$1$ on $X \otimes_k k_v$ also defines a retraction of this map,
by the formula $\mylangle A_v, \sum n_P P\myrangle_v = \sum n_P \Cores_{k_v(P)/k_v} A_v(P)$ for $A_v \in \Br(X \otimes_k k_v)$.
According to Lichtenbaum--Tate
duality~\cite[\textsection5]{lichtenbaum} for finite~$v$ and to a theorem of Witt~\cite[II', p.~5]{witt} for real~$v$ (see also~\cite[20.1.3]{scheiderer}),
all retractions of $\Br(k_v) \rightarrow \Br(X \otimes_k k_v)$
come from divisors of degree~$1$ in this way.
Hence for every place~$v$ of~$k$, there exists a degree~$1$ divisor $z_v$
on $X \otimes_k k_v$ such that $\mylangle A_v, z_v \myrangle_v=r_v(A_v)$ for all $A_v \in \Br(X \otimes_k k_v)$.

Let $A \in \Br(X)$.  Let~$\Omega$ denote the set of places of~$k$,
and for every $v \in \Omega$,  let $\inv_v \colon \Br(k_v) \hookrightarrow \Q/\Z$ denote the
canonical inclusion.  A glance at the construction of~$r$ and~$r_v$ reveals that the image of $r(A)$
in $\Br(k_v)$ coincides with $r_v(A \otimes_k k_v)$ for all $v\in\Omega$.  Since $r_v(A \otimes_k k_v)= \mylangle A\otimes_k k_v,z_v\myrangle_v$, it
follows that $\inv_v \mylangle A\otimes_k k_v,z_v\myrangle_v=\inv_v r(A)$.
Now $\sum_{v \in \Omega} \inv_v r(A)$ vanishes by the global reciprocity law;
hence
\begin{equation}
\label{suminvzero}
\sum_{v \in \Omega} \inv_v \mylangle A\otimes_k k_v,z_v\myrangle_v=0\text{.}
\end{equation}

Under the assumption that the Tate--Shafarevich group of the Jacobian of~$X$ is finite,
the existence of a family $(z_v)_{v \in \Omega}$ of local divisors of degree~$1$ satisfying~(\ref{suminvzero})
for every $A \in \Br(X)$ implies the existence of a divisor of degree~$1$ on the curve~$X$ itself, by a theorem
of Saito~\cite{saito} (see \cite[Prop.~3.7]{ctpspm}).  Thus the theorem is proved.
\end{proof}

\begin{rmks}
(i) Under the hypotheses of Theorem~\ref{thmab}, if one assumes that not
only~(\ref{seg}) but also~(\ref{segnonab}) splits, then there are many ways
to deduce that~$X$ admits a divisor of degree~$1$.
For instance, if~$s$ is a splitting of~(\ref{segnonab}) and
$\alpha \in H^2_\et(X,\Zhat(1))$ denotes its cycle class (see~\cite[Th.~2.6 and Rem.~2.7]{ew}),
the argument employed in the proof of~\cite[Prop.~3.1]{ew} shows that the image~$\beta$ of~$\alpha$
in the total Tate module $T(\Br(X))$ vanishes in $T(\Br(X \otimes_k k_v))$ for every place~$v$
of~$k$.
Now the finiteness assumption in Theorem~\ref{thmab}
ensures that $T(\Br(X))$ injects into $\prod_{v \in \Omega} T(\Br(X \otimes_k k_v))$,
so that $\beta=0$ and therefore~$\alpha$
belongs to the image of the cycle class map $\Pic(X) \otimes \Zhat \rightarrow H^2_\et(X,\Zhat(1))$;
which implies, as in~\cite[Cor.~3.6]{ew}, that~$X$ contains a divisor of degree~$1$.
The point of Theorem~\ref{thmab} is that the Galois-theoretic condition
which appears in its statement cannot be weakened since it is also a necessary condition
for the existence of divisors of degree~$1$.

(ii) Suppose~$k$ is a $p$\nobreakdash-adic field.  Using Roquette's theorem (see~\cite{lichtenbaum}),
it is easy to check that the statement of Theorem~\ref{thmab} still holds in this context.
However, more is true:
since,
by Lichtenbaum--Tate duality,
retractions of $\Br(k)\rightarrow \Br(X)$
are in one-to-one correspondence with
classes of degree~$1$ in $\Pic(X)$,
the proof of Theorem~\ref{thmab} shows that any splitting of~(\ref{seg}) determines
a well-defined class of degree~$1$ in $\Pic(X)$.
This is to be compared with Koenigsmann's theorem~\cite{koenigsmann} (recently refined by Pop~\cite{pop}
to a metabelian statement)
according to which any splitting of~(\ref{segnonab}) determines a well-defined
rational point on~$X$.
Of course, due to the abelian nature of~(\ref{seg}),
one cannot hope to associate rational points to splittings of~(\ref{seg}),
since there are curves over $p$\nobreakdash-adic fields which admit divisors of degree~$1$
but no rational points.

(iii) We do not know whether the statement of Theorem~\ref{thmab} still holds if~$k$ is only assumed to be a finitely generated field extension of~$\Q$.
For such a field~$k$, Grothendieck's section conjecture still predicts that a splitting of~(\ref{segnonab}) should determine a rational point of~$X$.
\end{rmks}

\section{Galois groups, fundamental groups and the elementary obstruction}
\label{sectob}

Let~$X$ now be an arbitrary smooth and geometrically irreducible variety over a field~$k$.  We keep the notation introduced
at the beginning of~\textsection\ref{sectionbirat}.  The exact sequence
\begin{equation}
\label{segbis}
\xymatrix{1 \ar[r] & \Gabgeom \ar[r] & \Gab \ar[r] & \Gk \ar[r] & 1}
\end{equation}
still makes sense, and it is still true that
it splits whenever~$X$ possesses a $0$\nobreakdash-cycle of degree~$1$, as a consequence of
the following lemma and a restriction-corestriction argument:

\begin{lem}
\label{deligne}
If $X(k)\neq\emptyset$ then the exact sequence~(\ref{introsegnonab}) splits.
\end{lem}

Lemma~\ref{deligne} is a straightforward generalisation of Deligne's remark for curves.  We include a proof for the convenience of the reader.

\begin{proof}
Let $K_0=k$ and for any $i \geq 1$, let $K_i=K_{i-1}((t_i))$, where the $t_i$'s are indeterminates.
For any $i\geq 1$, the natural projection of absolute Galois groups $G_{K_i} \to G_{K_{i-1}}$ admits a section (see \cite[Ch.~II, \textsection4.3, Ex.~1 and~2]{serrecg}).
Let $x \in X(k)$.
Since~$X$ is regular, the completion of the local ring of~$X$ at~$x$ is $k$\nobreakdash-isomorphic to $k[[t_1,\dots,t_n]]$.
It follows that~$k(X)$ embeds $k$\nobreakdash-linearly into~$K_n$.
By composing sections of $G_{K_i} \to G_{K_{i-1}}$ for $i \in \{1, \mskip2.5mu \ldots \mskip-1mu, \mskip2.5mu n\}$ with the projection $G_{K_n} \to G_{k(X)}$ given by such an embedding,
one obtains a splitting of~(\ref{introsegnonab}).
\end{proof}

According to Colliot-Th\'el\`ene and Sansuc~\cite[Prop.~2.2.2]{ctsandesc2}, the exact sequence of discrete $G_k$\nobreakdash-modules
\begin{equation}
\label{eqobelem}
\xymatrix{
0 \ar[r] & \kbarstar \ar[r] & \kbar(X)^\star \ar[r] & \kbar(X)^\star/\kbarstar \ar[r] & 0
}
\end{equation}
also splits ($G_k$\nobreakdash-equivariantly) whenever~$X$ possesses a $0$\nobreakdash-cycle of degree~$1$.
Thus the exact sequences~(\ref{segbis}) and~(\ref{eqobelem}) both constitute
obstructions to the existence of $0$\nobreakdash-cycles of degree~$1$ on~$X$.
The former appears in the statement of Theorem~\ref{thmab} while the latter
visibly plays a r\^ole in the proof of the same theorem (see~(\ref{diagrammem})),
which prompts the question of their exact relation.
We address this question in Theorem~\ref{thmcomp} below.

Let $\ob(X) \in \Ext^1_{G_k}(\kbar(X)^\star/\kbarstar,\kbarstar)$ denote the class of the extension~(\ref{eqobelem}).
Following the authors of~\cite{ctsandesc2}
we call $\ob(X)$ the \emph{elementary obstruction} (to the existence of $0$\nobreakdash-cycles of degree~$1$ on~$X$).

\begin{thm}
\label{thmcomp}
Let~$X$ be a smooth and geometrically irreducible variety over a field~$k$ of characteristic~$0$.
The exact sequence of profinite groups~(\ref{segbis}) splits if and only if $\ob(X)$ belongs
to the maximal divisible subgroup
of $\Ext^1_{G_k}(\kbar(X)^\star/\kbarstar,\kbarstar)$.
\end{thm}

Under the additional assumption that $\Pic(X \otimes_k \kbar)$ is finitely
generated, which notably excludes curves of positive genus,
the conclusion of Theorem~\ref{thmcomp} was essentially known (see~\cite[Cor.~2.4.2]{skobook} and Lemma~\ref{picfg} below).

\begin{proof}
To establish Theorem~\ref{thmcomp}, we shall exploit the
commutative diagram of discrete $\Gab$\nobreakdash-modules
already encountered in the proof of Theorem~\ref{thmab}, namely
\begin{equation}
\begin{aligned}
\label{diagrammembis}
\xymatrix{
0 \ar[r] & \kbarstar \ar@{=}[d] \ar[r] & \kbar(X)^\star \ar[d]^i \ar[r] & \ar[d]\kbar(X)^\star/\kbarstar \ar[r] & 0 \\
0 \ar[r] & \kbarstar \ar[r] & M \ar[r]^(0.3)p & \left(\kbar(X)^\star/\kbarstar\right) \otimes_\Z \Q \ar[r] & 0
}
\end{aligned}
\end{equation}
where
$M = \bigset{f \in \kbar(X)^{\ab \mkern1mu \star}}{\exists\mkern1mu n \geq 1, \mkern2mu{}f^n \in \kbar(X)^\star}$.
One half of the statement of Theorem~\ref{thmcomp} is a more or less immediate
consequence of the existence of this diagram.
If~(\ref{segbis}) splits, then, as in the proof of
Theorem~\ref{thmab}, we may consider~(\ref{diagrammembis}) as a commutative
diagram of discrete $G_k$\nobreakdash-modules.
The top row of~(\ref{diagrammembis}) is then the pullback of the bottom row of~(\ref{diagrammembis})
by the natural map $\kbar(X)^\star/\kbarstar \rightarrow \left(\kbar(X)^\star/\kbarstar\right) \otimes_\Z \Q$
not only in the category of abelian groups but also in the category of discrete $G_k$\nobreakdash-modules.
In other words $\ob(X)$ is the image, by pullback, of
the class of the bottom row of~(\ref{diagrammembis}) in
$\Ext^1_{G_k}((\kbar(X)^\star/\kbarstar) \otimes_\Z \Q,\kbarstar)$.  The latter group being a $\Q$\nobreakdash-vector
space, it follows that $\ob(X)$ belongs to the maximal divisible subgroup of
 $\Ext^1_{G_k}(\kbar(X)^\star/\kbarstar,\kbarstar)$.

For the converse implication we start with a few lemmas.

\begin{lem}
\label{lemsec}
Splittings of the exact sequence of profinite groups~(\ref{segbis}) are canonically in one-to-one correspondence with discrete
$G_k$\nobreakdash-module structures on the abelian group~$M$ which extend the natural $G_k$\nobreakdash-module structure
of the subgroup $\kbar(X)^\star$ of~$M$.
\end{lem}

\begin{proof}
A splitting of~(\ref{segbis}) induces an action of~$G_k$ on~$\kbar(X)^{\ab \mkern1mu \star}$, hence on~$M$, satisfying the required property.  Suppose
conversely that~$M$ is endowed with such a structure.

We claim that every
 $\sigma \in G_k$
admits a unique lifting
$\tilde\sigma \in \smash{\Gab}$ such that the given action of~$\sigma$ on~$M$ coincides with the natural
action of~$\tilde\sigma$ on $M \subset \kbar(X)^{\ab \mkern1mu \star}$.
Assuming this, the map $\Gk \rightarrow \Gab$ which sends~$\sigma$ to $\tilde\sigma$
is a continuous group homomorphism and therefore constitutes a splitting of~(\ref{segbis}).  We are thus reduced to proving the claim.

Let $\sigma \in G_k$.  Choose an arbitrary lifting~$\tilde\sigma$ of~$\sigma$.
The endomorphism $\varphi \colon M \rightarrow M$
defined by $\varphi(f)=\tilde\sigma^{-1}(\sigma(f))/f$
vanishes on $\kbar(X)^\star$ and takes values in the torsion subgroup~$\mmu_\infty$ of~$\kbarstar$, so that it factors as
$$
M \rightarrow M/\kbar(X)^\star \xrightarrow{\;\bar \varphi}  \mmu_\infty \subset \kbarstar \subset M\text{.}
$$
Now by Kummer theory, one has
$$M/\kbar(X)^\star=(\kbar(X)^\star/\kbarstar) \otimes_\Z \Q/\Z=\Hom(\Gabgeom,\mmu_\infty)$$
(where~$\Hom$ denotes the set of continuous homomorphisms with respect to the discrete topology on $\mmu_\infty$).
Hence~$\bar \varphi$ may be seen as a map $\Hom(\Gabgeom,\mmu_\infty) \rightarrow \mmu_\infty$.
By Pontrjagin duality we conclude that there is a unique $\tau \in \Gabgeom$ such that $\bar\varphi(f)=\tau(f)/f$ for all $f \in M$.
The sought for lifting of~$\sigma$ is then $\tilde\sigma\tau$ and it is indeed unique since~$\varphi$ determines~$\tau$.
\end{proof}

\begin{lem}
\label{lemsecconj}
Splittings of~(\ref{segbis}) up to conjugation
by~$\Gabgeom$ are canonically in one-to-one correspondence with
triples $(E,i_E,p_E)$ up to isomorphism, where~$E$
is a discrete $G_k$\nobreakdash-module and~$i_E$ and~$p_E$ are $G_k$\nobreakdash-equivariant
maps fitting into a commutative diagram
\begin{equation}
\begin{aligned}
\label{diagrammee}
\xymatrix{
0 \ar[r] & \kbarstar \ar@{=}[d] \ar[r] & \kbar(X)^\star \ar[d]^{i_E} \ar[r] & \ar[d]\kbar(X)^\star/\kbarstar \ar[r] & 0 \\
0 \ar[r] & \kbarstar \ar[r] & E \ar[r]^(0.3){p_E} & \left(\kbar(X)^\star/\kbarstar\right) \otimes_\Z \Q \ar[r] & 0
}
\end{aligned}
\end{equation}
with exact rows. (By an isomorphism of triples $(E,i_E,p_E) \isoto (E',i_{E'},p_{E'})$, we mean an isomorphism $h \colon E \isoto E'$ such that
$h \circ i_E=i_{E'}$ and $p_{E'}\circ h = p_E$.)
\end{lem}

\begin{proof}
Lemma~\ref{lemsecconj} is a formal consequence of Lemma~\ref{lemsec} once one has verified
that for any triple $(E,i_E,p_E)$, the diagram~(\ref{diagrammee}) is isomorphic, as a diagram of abelian groups
(disregarding the action of~$G_k$),
to~(\ref{diagrammembis}); indeed the choice of such an isomorphism allows one to transport the $G_k$\nobreakdash-module structure from~$E$
to~$M$.

To this end, we first note that the exactness of the second row of~(\ref{diagrammee}) 
forces~$E$ to be divisible, so that there exists a homomorphism $h \colon M \rightarrow E$ satisfying $h \circ i = i_E$.
The kernel and cokernel of~$h$ are simultaneously $\Q$\nobreakdash-vector spaces (since
$E/\kbarstar$ and $M/\kbarstar$ are $\Q$\nobreakdash-vector spaces)
and torsion groups (since $E/\kbar(X)^\star$ and
$M/\kbar(X)^\star$ are torsion). Hence~$h$ is an isomorphism.
Moreover, we automatically have $p=p_E\circ h$
since the homomorphism $p-(p_E\circ h)$ vanishes on the image of~$i$ and takes values in a torsionfree group, and the cokernel of~$i$ is torsion.
\end{proof}

\begin{lem}
\label{lemprofini}
Let~$G$ be a profinite group and~$A$, $B$ be two discrete $G$\nobreakdash-modules.
The image of the natural map
$\Ext^1_G(A \otimes_\Z \Q,B) \rightarrow \Ext^1_G(A,B)$
is the maximal divisible subgroup of $\Ext^1_G(A,B)$.
\end{lem}

\begin{proof}
Since $\Ext^1_G(A \otimes_\Z \Q,B)$ is a $\Q$\nobreakdash-vector space,
this map naturally factors as
$$
\xymatrix{
\Ext^1_G(A \otimes_\Z\Q,B) \ar[r]^(0.45)u & \Hom(\Q,\Ext^1_G(A,B)) \ar[r]^(0.6)v & \Ext^1_G(A,B)\text{,}
}
$$
where~$v$ is the evaluation at~$1$ map.  The image of~$v$ is the maximal divisible subgroup of $\Ext^1_G(A,B)$;
hence it suffices to check that~$u$ is onto.  In view of the spectral sequence
$\Ext^p(\Q,\Ext^q_G(A,B))\Rightarrow \Ext^{p+q}_G(A \otimes_\Z \Q,B)$ (see~\cite[Ch.~I, \textsection0]{milneadt}),
the cokernel of~$u$ embeds into $\Ext^2(\Q,\Hom_G(A,B))$. But $\Ext^2(\Q,\Hom_G(A,B))=0$ as the category of abelian groups has global dimension~$1$
(see, \emph{e.g.}, \cite[Ch.~VI, Prop.~2.8]{cartaneilenberg}).
\end{proof}

We are now in a position to complete the proof of Theorem~\ref{thmcomp}.
Assume $\ob(X)$ belongs to the maximal divisible subgroup of
$\Ext^1_{G_k}(\kbar(X)^\star/\kbarstar,\kbarstar)$.  According to Lemma~\ref{lemprofini}, there exists a discrete $G_k$\nobreakdash-module~$E$
and $G_k$\nobreakdash-equivariant maps~$i_E$ and~$p_E$ such that the diagram~(\ref{diagrammee}) commutes and has exact rows.
Lemma~\ref{lemsecconj} now implies that~(\ref{segbis}) splits.
\end{proof}

\begin{rmk}
\label{imperfect}
Let~$k$ be a field of characteristic $p>0$. For any prime number~$\ell$,
denote by $\Gabgeomell$ the $\ell$\nobreakdash-Sylow subgroup of~$\Gabgeom$
and by~$D$ (resp.~$D_\ell$) the maximal divisible (resp.~$\ell$\nobreakdash-divisible) subgroup
of $\Ext^1_{G_k}(\kbar(X)^\star/\kbarstar,\kbarstar)$.
Let $\Gabgeomprimetop = \prod_{\ell \neq p}\Gabgeomell$ and $D_{p'}=\bigcap_{\ell \neq p}D_\ell$,
so that $\Gabgeom=\Gabgeomp \oplus \Gabgeomprimetop$ and $D=D_p \cap D_{p'}$.
In this situation, the proof of Theorem~\ref{thmcomp} given above is easily adapted to show that
the exact sequence obtained by pushing out~(\ref{segbis}) by the projection $\Gabgeom \to \Gabgeomprimetop$
splits if and only if $\ob(X) \in D_{p'}$.
On the other hand, the exact sequence obtained by pushing out~(\ref{segbis}) by the projection $\Gabgeom \to \Gabgeomp$ always splits
(see~\cite[Ch.~II, \textsection2.2, Prop.~3 and Ch.~I, \textsection3.4, Prop.~16]{serrecg}).
Moreover, if~$k$ is perfect, multiplication by~$p$
is an automorphism of the abelian group $\Ext^1_{G_k}(\kbar(X)^\star/\kbarstar,\kbarstar)$
since it is an automorphism of the $G_k$\nobreakdash-module~$\kbarstar$. Therefore
in this case $D=D_{p'}$.
From these remarks we deduce that the statement of Theorem~\ref{thmcomp} holds, more generally, over perfect fields of arbitrary characteristic.

Over imperfect fields, however, the conclusion of Theorem~\ref{thmcomp} may fail.
Indeed, let~$k$ be a field of cohomological dimension~$1$
which is not of dimension~$\leq 1$ in the sense of Serre~\cite[Ch.~II, \textsection3]{serrecg}.  Then~(\ref{segbis})
splits for any~$X$ (by~\cite[Ch.~I, \textsection5.9, Cor.~2]{serrecg}) whereas there exists a smooth and geometrically irreducible toric variety~$X$ over~$k$ such that $\ob(X)\neq 0$
(see~\cite[Prop.~3.4.3]{w}).  Since~$X$ is toric, the following lemma shows that $\ob(X)$ does not belong to~$D$.
\end{rmk}

\begin{lem}
\label{picfg}
Let~$X$ be a smooth and geometrically irreducible variety over a field~$k$.
If $\Pic(X \otimes_k \kbar)$ is finitely generated then
$\Ext^1_{G_k}(\kbar(X)^\star/\kbarstar,\kbarstar)$ has finite exponent.
\end{lem}

\begin{proof}
After shrinking~$X$ we may assume that $\Pic(X \otimes_k \kbar)=0$. Denote by $\kbar[X]^\star$ the group
of regular invertible functions on $X\otimes_k \kbar$
and by~$T$ the torus over~$k$ with character group $\kbar[X]^\star/\kbarstar$.
The group $\Ext^1_{G_k}(\kbar[X]^\star/\kbarstar,\kbarstar)$
is isomorphic to $H^1(k,T)$ (see \cite[Lemma~2.3.7]{skobook}), so that
by Hilbert's Theorem~90 it is killed by the degree of any finite extension of~$k$ which splits~$T$.
On the other hand $\Ext^1_{G_k}(\Div(X \otimes_k \kbar),\kbarstar)$ vanishes by Shapiro's lemma and Hilbert's Theorem~90,
since $\Div(X \otimes_k \kbar)$ is a permutation $G_k$\nobreakdash-module.
In view of the exact sequence
$$
\xymatrix{
0 \ar[r] & \kbar[X]^\star/\kbarstar \ar[r] & \kbar(X)^\star/\kbarstar \ar[r] & \Div(X \otimes_k \kbar) \ar[r] & 0 \rlap{\text{,}}
}
$$
the lemma follows.
\end{proof}

From Theorem~\ref{thmcomp} and Remark~\ref{imperfect} one deduces:

\begin{cor}
\label{corsplit}
Let~$X$ be a smooth and geometrically irreducible variety over a field~$k$.
If the elementary obstruction to the existence of a $0$\nobreakdash-cycle of degree~$1$ on~$X$
vanishes, then the exact sequence of profinite groups~(\ref{segbis}) splits.
\end{cor}

The statement that if~$X$ is a proper variety over a field of characteristic~$0$ and the elementary obstruction vanishes
then the abelianisation of~(\ref{introsepi}) splits has independently been
shown to hold by Harari and Szamuely~\cite[Rem.~5.6]{hasza}.
It also follows from Corollary~\ref{corsplit} since splittings of~(\ref{segbis})
induce splittings of the abelianisation of~(\ref{introsepi}).

Theorem~\ref{thmcomp} is sharp in the sense that it characterises
the existence of ``birational abelian sections'' (that is, of splittings of~(\ref{segbis}))
purely in terms of the elementary obstruction.
Obviously the condition of Theorem~\ref{thmcomp} is satisfied if the elementary obstruction vanishes;
however the question remains whether there may exist fields~$k$ and smooth and geometrically irreducible varieties~$X$ such that
$\ob(X)$ belongs to the maximal divisible subgroup
of $\Ext^1_{G_k}(\kbar(X)^\star/\kbarstar,\kbarstar)$ without vanishing.  We do not know of any such example.
Neither the field~$k$ nor the variety~$X$ can be too simple for this to occur.  Indeed
the abelian group $\Pic(X \otimes_k \kbar)$ cannot be finitely generated
(by Lemma~\ref{picfg}),
the field~$k$ cannot have dimension~$\leq 1$ (according to~\cite[Th.~3.4.1]{w}), and in Theorem~\ref{thmobobdiv} below we prove that~$k$ cannot be a $p$\nobreakdash-adic
field, a real closed field, or a number field, assuming in the latter case the finiteness of a Tate--Shafarevich group.

\begin{thm}
\label{thmobobdiv}
Let~$X$ be a smooth, proper and geometrically irreducible variety over a field~$k$.  Assume either that~$k$ is a real closed field or a $p$\nobreakdash-adic field,
or that~$k$ is a number field and the Tate--Shafarevich group of the Picard variety of~$X$ is finite.
Then the exact sequence of profinite groups~(\ref{segbis}) splits if and only if $\ob(X)=0$.
\end{thm}

The argument we give for Theorem~\ref{thmobobdiv} in the case of number fields was to a large extent inspired by the proof of \cite[Th.~2.13]{boctsk}.
It also relies on Theorem~\ref{thmcomp}, on the strategy used in the proof of Theorem~\ref{thmab}, on~\cite[Th.~1.1]{haszlocglob}
and on \cite[Th.~2.2]{w}.

\begin{proof}
If $\ob(X)=0$, then~(\ref{segbis}) splits by Theorem~\ref{thmcomp}. Conversely, suppose~(\ref{segbis}) splits and let $s \colon \Gk \rightarrow \Gab$ be a splitting.
If~$V$ is a geometrically integral variety over a field~$K$, we set $\Br_1(V)=\Ker\big(\Br(V) \to \Br(V \otimes_K \Kbar)\big)$,
where~$\Kbar$ denotes an algebraic closure of$~K$,
and $\Br_1(K(V))=\Ker\big(\Br(K(V)) \to \Br(\Kbar(V))\big)$.
There is a canonical isomorphism $\Br_1(k(X)) = H^2(k,\kbar(X)^\star)$,
so that as in the proof of Theorem~\ref{thmab}, the section~$s$ induces a retraction $r \colon \Br_1(k(X)) \rightarrow \Br(k)$ of the natural map $\Br(k) \rightarrow \Br_1(k(X))$.
The latter is therefore injective, which, if~$k$ is a real closed field or a $p$\nobreakdash-adic field, implies that $\ob(X)=0$ (see~\cite[Th.~2.5 and~2.6]{boctsk}).
Thus we may assume that~$k$ is a number field.  Then, again as in the proof of Theorem~\ref{thmab},
for any place~$v$ of~$k$
the section~$s$ induces a retraction $r_v \colon \Br_1(k_v^h(X)) \rightarrow \Br(k_v^h)$
of the natural map $\Br(k_v^h) \rightarrow \Br_1(k_v^h(X))$.
It satisfies $\inv_v r_v(A \otimes_k k_v^h) = \inv_v r(A)$ for all $A \in \Br_1(k(X))$, where $\inv_v \colon \Br(k_v^h) \hookrightarrow \Q/\Z$ denotes the invariant map
given by local class field theory.
Let $U \subseteq X$ be a dense open subset and let~$\Alb^1_{U/k}$ denote the Albanese torsor of~$U$ (see \cite[\textsection2]{w}; it is a torsor under a semi-abelian variety over~$k$).
This variety comes by definition with a morphism $u \colon U \rightarrow \Alb^1_{U/k}$;
as a consequence~$r$ and~$r_v$, for any place~$v$, induce retractions
$r' \colon \Br_1(\Alb^1_{U/k}) \rightarrow \Br(k)$ and $r'_v \colon \Br_1(\Alb^1_{U/k} \otimes_k k_v^h) \rightarrow \Br(k_v^h)$ of the natural maps
in the other direction.
  In particular $\Br(k_v^h)$ injects into $\Br(\Alb^1_{U/k} \otimes_k k_v^h)$ for every place~$v$ of~$k$
and therefore $\Alb^1_{U/k}$ possesses an adelic point
(see~\cite[Th.~2.5, Th.~2.6 and Th.~3.2]{boctsk}).  Fix an adelic point  $(P_v)_{v \in \Omega}$ of $\Alb^1_{U/k}$.
Denote by~$\Ba(\Alb^1_{U/k})$ the kernel of $\Br_1(\Alb^1_{U/k}) \rightarrow \prod_{v \in \Omega} \Br(\Alb^1_{U/k}\otimes_k k_v)/\Br(k_v)$.
We have
 $\Ba(\Alb^1_{U/k})=\Ker\big(\Br_1(\Alb^1_{U/k}) \rightarrow \prod_{v \in \Omega} \Br(\Alb^1_{U/k}\otimes_k k_v^h)/\Br(k_v^h)\big)$
since $\Br(k_v^h)=\Br(k_v)$
and the natural map $\Br(X \otimes_k k_v^h) \rightarrow \Br(X \otimes_k k_v)$ is injective
(see~\cite[p.~334]{boctsk}).
  As a consequence, and in view of the fact that~$r'_v$ is a retraction of $\Br(k_v^h) \rightarrow \Br_1(\Alb^1_{U/k}\otimes_k k_v^h)$,
we see that $\inv_v A(P_v)=\inv_v r'_v(A \otimes_k k_v^h)$ for all $A \in \Ba(\Alb^1_{U/k})$ and all~$v$.
On the other hand, $\inv_v r'_v(A\otimes_k k_v^h)=\inv_v r'(A)$.  Hence, by global reciprocity, we find that $\sum_{v \in \Omega}\inv_v A(P_v)=0$.
In other words, the adelic point $(P_v)_{v \in \Omega}$ of~$\Alb^1_{U/k}$ is orthogonal to $\Ba(\Alb^1_{U/k})$ with respect to the Brauer--Manin pairing.
By~\cite[Th.~1.1]{haszlocglob}, it then follows that $\Alb^1_{U/k}(k)\neq\emptyset$.
Since~$U$ was arbitrary, we may now apply \cite[Th.~2.2]{w} and conclude that $\ob(X)=0$.
\end{proof}

We conclude this note by addressing a question about the elementary obstruction which was raised by Borovoi, Colliot-Th\'el\`ene and Skorobogatov in~\cite[p.~327]{boctsk} and
which has no apparent connection with Galois groups and fundamental groups.

\begin{question}
\label{qbctsk}
Let~$X$ be a geometrically integral variety over a field~$k$.  Let~$K/k$ be a field extension.  Does $\ob(X)=0$ imply $\ob(X \otimes_k K)=0$~?
\end{question}

Partial positive answers were given in \cite[Prop.~2.3]{boctsk} and \cite[Cor.~3.2.3, Cor.~3.3.2, Prop.~3.4.4]{w}.
As it turns out, considering a situation in which the analogue of Grothendieck's section conjecture
thoroughly fails leads one to a negative answer to Question~\ref{qbctsk}:

\begin{thm}
\label{thmcex}
Let $k=\C((t))$.  There exists a geometrically integral curve~$X$ over~$k$ and a field extension $K/k$
such that $\ob(X)=0$ but $\ob(X \otimes_k K) \neq 0$.
\end{thm}

\begin{proof}
Let~$C$ be the projective plane curve over~$k$ defined by $x^3+ty^3+t^2z^3=0$.
Obviously~$C$ has no rational points.
Since~$C$ is a smooth and proper curve of genus~$1$, it follows
(by the Riemann--Roch theorem) that there is no divisor of degree~$1$ on~$C$.
On the other hand~$C$ has closed points of degree~$3$.
We conclude that the degree of any divisor on~$C$ must be divisible by~$3$.

The curve~$C$ being geometrically connected, its fundamental group $\pi_1(C)$ surjects onto the absolute Galois group~$G_k$
of~$k$, which is isomorphic to~$\Zhat$.  The choice of a lifting of $1 \in \Zhat$ to $\pi_1(C)$ determines a section
$s \colon G_k \rightarrow \pi_1(C)$ of the projection $\pi_1(C) \rightarrow G_k$.
Let $\alpha \in H^2_\et(C,\mmu_3)$ denote the cycle class of~$s$ in \'etale cohomology with $\Z/3\Z$ coefficients
(see~\cite[Th.~2.6]{ew}).
We recall from \cite[Th.~2.6]{ew} that~$\alpha$ has degree~$1$ modulo~$3$ (in the sense that its image
in $H^2_\et(C \otimes_k \kbar,\mmu_3)=\Z/3\Z$ is equal to~$1$).
Therefore~$\alpha$ cannot be the cycle class of a divisor on~$C$.

We also recall that by the construction of~$\alpha$, there exists an \'etale cover $\pi \colon X \rightarrow C$
such that~$X$ is geometrically irreducible over~$k$ and
such that the class $p^\star \alpha \in H^2_\et(X \times C,\mmu_3)$, where $p \colon X \times C \rightarrow C$ denotes the second projection,
is equal to the cycle class of the graph of~$\pi$.  

The Kummer exact sequence
\smash{$\xymatrix@C=4.2ex{1 \ar[r] & \mmu_3 \ar[r] & \Gm \ar[r]\ar@{}[r]<-.15em>^(.47){\times 3} & \Gm \ar[r] & 1}$}
gives rise to a commutative diagram with exact rows
\begin{equation}
\label{cexcommdiag}
\begin{aligned}
\xymatrix{
 \Pic(C) \ar[d] \ar[r] & H^2_\et(C,\mmu_3) \ar[d]\ar[r] & \Br(C) \ar[d] \\
 \Pic(X \times C) \ar[r] & H^2_\et(X\times C,\mmu_3) \ar[r] & \Br(X\times C)\rlap{\text{,}}
}
\end{aligned}
\end{equation}
where the vertical arrows are the pullback maps~$p^\star$.
Since the cycle class of the graph of~$\pi$ in $H^2_\et(X\times C,\mmu_3)$ comes from $\Pic(X \times C)$, the above diagram shows that the image of~$\alpha$ in $\Br(X \times C)$ vanishes.
On the other hand the image~$\beta$ of~$\alpha$ in $\Br(C)$ does not vanish, since~$\alpha$ is not the cycle class of a divisor.
Hence $\beta$ is a nonzero element of the kernel of $p^\star \colon \Br(C) \rightarrow \Br(X \times C)$.
As a consequence, this map is not injective, and
the natural map $\Br(K) \rightarrow \Br(X \otimes_k K)$, where $K=k(C)$, is not injective either.
In particular $\ob(X \otimes_k K) \neq 0$ (see~\cite[Prop.~2.2.5]{ctsandesc2}).
Finally we have $\ob(X)=0$ because~$k$ is a field of dimension~$\leq 1$ (see~\cite[Th.~3.4.1]{w}).
\end{proof}

\begin{rmks}
\label{rmkfin}
(i) An alternative way of deducing Theorem~\ref{thmcex} from~\cite[Th.~3.4.1]{w} is given in \cite[Prop.~2.4.2]{cipkrash}.
We are grateful to Daniel Krashen for this remark.

(ii) In the example of Theorem~\ref{thmcex}, not only does the elementary obstruction on $X \otimes_k K$ not vanish;
the class $\ob(X \otimes_k K)$ does not even belong to the maximal divisible subgroup of
$\Ext^1_{G_K}(\Kbar(X)^\star/\Kbarstar,\Kbarstar)$.

Indeed, if it were the case, then by Theorem~\ref{thmcomp},
the natural projection $\GabK \rightarrow \GK$ would admit a section. Therefore, by the argument used at the beginning of the proof of Theorem~\ref{thmab},
the map $\Br(K) \rightarrow \Br(X \otimes_k K)$ would have to be injective, which, as we have seen, it is not.

(iii) According to Remark~\ref{rmkfin}~(ii) and to Theorem~\ref{thmcomp}, the proof of Theorem~\ref{thmcex} also gives an example
of a geometrically integral curve~$X$ over a field~$k$ of characteristic~$0$, and of a field extension~$K/k$,
such that the exact sequences~(\ref{introsegnonab}) and~(\ref{introsegab}) associated to the curve~$X$ over the field~$k$
both split, but such that neither of the exact sequences~(\ref{introsegnonab}) and~(\ref{introsegab}) associated to the curve $X \otimes_k K$ over the field~$K$
splits.
\end{rmks}

\section*{Acknowledgments}

Nguy\~{\^e}n Duy T\^an pointed out an inaccuracy in the proof
of Theorem~\ref{thmobobdiv} in an earlier version of this note.
Jakob Stix reminded us that fields of characteristic $p>0$
have $p$\nobreakdash-cohomological dimension~$1$, which was the starting
point for Remark~\ref{imperfect}.  We are grateful to both of them.

\bibliographystyle{amsplain}
\bibliography{bir_ab}

\providecommand{\bysame}{\leavevmode\hbox to3em{\hrulefill}\thinspace}
\providecommand{\MR}{\relax\ifhmode\unskip\space\fi MR }
% \MRhref is called by the amsart/book/proc definition of \MR.
\providecommand{\MRhref}[2]{%
  \href{http://www.ams.org/mathscinet-getitem?mr=#1}{#2}
}
\providecommand{\href}[2]{#2}
\begin{thebibliography}{10}

\bibitem{boctsk}
M.~Borovoi, J.-L. Colliot-Th{\'e}l{\`e}ne, and A.~N. Skorobogatov, \emph{The
  elementary obstruction and homogeneous spaces}, Duke Math. J. \textbf{141}
  (2008), no.~2, 321--364.

\bibitem{cartaneilenberg}
H.~Cartan and S.~Eilenberg, \emph{Homological algebra}, Princeton University
  Press, Princeton, NJ, 1956.

\bibitem{cipkrash}
M.~\c{C}iperiani and D.~Krashen, \emph{Relative {B}rauer groups of genus~{$1$}
  curves}, preprint, 2007.

\bibitem{ctpspm}
J-L. Colliot-Th{\'e}l{\`e}ne, \emph{Conjectures de type local-global sur
  l'image des groupes de {C}how dans la cohomologie \'etale}, Algebraic
  {$K$}-theory ({S}eattle, {WA}, 1997), Proc. Sympos. Pure Math., vol.~67,
  Amer. Math. Soc., Providence, RI, 1999, pp.~1--12.

\bibitem{ctsandesc2}
J-L. Colliot-Th{\'e}l{\`e}ne and J-J. Sansuc, \emph{La descente sur les
  vari\'et\'es rationnelles. {II}}, Duke Math. J. \textbf{54} (1987), no.~2,
  375--492.

\bibitem{delignegroupefond}
P.~Deligne, \emph{Le groupe fondamental de la droite projective moins trois
  points}, Galois groups over {${\bf Q}$} ({B}erkeley, {CA}, 1987), MSRI Publ.,
  vol.~16, Springer, New York, 1989, pp.~79--297.

\bibitem{ew}
H.~Esnault and O.~Wittenberg, \emph{Remarks on cycle classes of sections of the
  arithmetic fundamental group}, Mosc. Math. J. \textbf{9} (2009), no.~3,
  451--467.

\bibitem{gtof}
A.~Grothendieck, letter to {F}altings dated {J}une 27, 1983 (in German),
  published in: Geometric Galois actions, Vol.~1, London Math. Soc. Lecture
  Note Ser., vol.~242, Cambridge Univ. Press, Cambridge, 1997, pp.~49--58.

\bibitem{grothbr2}
\bysame, \emph{Le groupe de {B}rauer {II} : th\'eorie cohomologique}, Dix
  expos\'es sur la cohomologie des sch\'emas, Advanced studies in pure
  mathematics, vol.~3, Masson \& North-Holland, Paris, Amsterdam, 1968,
  pp.~67--87.

\bibitem{haszlocglob}
D.~Harari and T.~Szamuely, \emph{Local-global principles for 1-motives}, Duke
  Math. J. \textbf{143} (2008), no.~3, 531--557.

\bibitem{hasza}
\bysame, \emph{Galois sections for abelianized fundamental groups}, Math. Ann.
  \textbf{344} (2009), no.~4, 779--800, with an appendix by E. V. Flynn.

\bibitem{koenigsmann}
J.~Koenigsmann, \emph{On the `section conjecture' in anabelian geometry}, J.
  reine angew. Math. \textbf{588} (2005), 221--235.

\bibitem{lichtenbaum}
S.~Lichtenbaum, \emph{Duality theorems for curves over {$p$}-adic fields},
  Invent. math. \textbf{7} (1969), 120--136.

\bibitem{milneadt}
J.~S. Milne, \emph{Arithmetic duality theorems}, Perspectives in Mathematics,
  vol.~1, Academic Press Inc., Boston, MA, 1986.

\bibitem{pop}
F.~Pop, \emph{On the birational {$p$}-adic section conjecture}, to appear in
  Compositio Math.

\bibitem{saito}
S.~Saito, \emph{Some observations on motivic cohomology of arithmetic schemes},
  Invent. math. \textbf{98} (1989), no.~2, 371--404.

\bibitem{scheiderer}
C.~Scheiderer, \emph{Real and \'etale cohomology}, Lecture Notes in
  Mathematics, vol. 1588, Springer-Verlag, Berlin, 1994.

\bibitem{serrecg}
J-P. Serre, \emph{Cohomologie galoisienne}, fifth ed., Lecture Notes in
  Mathematics, vol.~5, Springer-Verlag, Berlin, 1994.

\bibitem{skobook}
A.~N. Skorobogatov, \emph{Torsors and rational points}, Cambridge Tracts in
  Mathematics, vol. 144, Cambridge University Press, Cambridge, 2001.

\bibitem{witt}
E.~Witt, \emph{{Z}erlegung reeller algebraischer {F}unktionen in {Q}uadrate.\
  {S}chiefk\"orper \"uber reellem {F}unktionenk\"orper}, J. reine angew. Math.
  \textbf{171} (1934), 4--11.

\bibitem{w}
O.~Wittenberg, \emph{On {A}lbanese torsors and the elementary obstruction},
  Math. Ann. \textbf{340} (2008), no.~4, 805--838.

\end{thebibliography}
\end{document}